# Ivan Boldyrev

(Radboud University, Nijmegen, The Netherlands)

# Soviet Mathematics and Economic Theory in the Past Century: An Historical Reappraisal[1]

What are the effects of authoritarian regimes on scholarly research in economics? And how might economic theory survive ideological pressures? The article addresses these questions by focusing on the mathematization of economics over the past century and drawing on the history of Soviet science. Mathematics in the USSR remained internationally competitive and generated many ideas that were taken up and played important roles in economic theory. These same ideas, however, were disregarded or adopted only in piecemeal fashion by Soviet economists, despite the efforts of influential scholars to change the economic research agenda. The article draws this contrast into sharper focus by exploring the work of Soviet mathematicians in optimization, game theory, and probability theory that was used in Western economics. While the intellectual exchange across the Iron Curtain did help advance the formal modeling apparatus, economics could only thrive in an intellectually open environment absent under the Soviet rule.

## 1. Introduction

While Soviet mathematics was on a par with the rest of the world, Soviet economics – especially after the Stalinist repressions of the 1930s – nearly disappeared as an academic discipline. This paper argues that exploring the adoption of Soviet mathematical work in economic theory helps illuminate how ideological intolerance and political oppression influenced economic thinking of the past century.

Historians of economics – unlike economic historians – rarely possess the data that would allow for a full-fledged counterfactual analysis. Of course, a historical exploration can be valuable in itself. In reconstructing interdisciplinary contexts and clarifying scholarly biographies, it adds new episodes to the broader history of ideas. In some cases, however, this exploration can also provide evidence for certain counterfactual claims.

The Soviet case is particularly salient in this respect. In the 1920s, many economists working in the USSR, just like mathematicians, communicated on equal terms with the international leaders of their fields. This was especially true for the members of the Moscow Conjuncture Institute (Klein 1999). The work of the institute, headed by Nikolai Kondratiev, arguably marks the heyday of early Soviet economics. However, the demise of the institute in 1928 and the repression of many of its members, including Kondratiev himself, were severe blows. Within a decade, most of those who were not repressed, like Eugen Slutsky, stopped doing economics (Barnett 2011). It is hard to find an example of an environment where the culture of teaching economics and doing economic research changed so drastically and where economics was subjected to such a long-term ideological pressure (Yasny 1972).

[1] I am deeply grateful to the editors of this journal, Steven Durlauf and David Romer, as well as to four anonymous referees for their comments that helped to significantly improve the paper. At different stages of this project, I learned a lot from the conversations with Rein van Alebeek, Valentina Desnitskaya, Mark Levin, Serge Lvovsky, Leonid Polishchuk, Victor Polterovich, Boris Polyak, Konstantin Sonin, the late Anatoly Vershik, and Elena Yanovskaya. The various support I received from Vladimír Dlouhý, Maryanne Johnson, Alexander Kouriaev, and Zachary Reyna is also gratefully acknowledged. Any errors or omissions remain my sole responsibility.

While economic research collapsed, Soviet mathematics thrived despite its post-1930s isolation. And while some individual mathematicians were targeted during Soviet times (Seneta 2004; Fuchs 2007; Graham and Kantor 2009; Hollings 2013), this did not suppress the field *in its entirety*. In short, mathematics did not suffer as economics did.[2] Moreover–albeit this is less well known–Soviet mathematics continued to produce ideas that were actively taken up in Western economic theory.

Formal techniques exerted an (often invisible) influence on the economics of the past century. Accounting for the mathematization of economics becomes particularly instructive once we detail how it happened in different ideological contexts. The history of mathematics thus serves as a foil to re-evaluate the differences in the ways economic ideas developed in the Soviet Union, a non-democracy with a command economy, and in the Western world. Documenting how the work of Soviet mathematicians was received in economic theory *beyond their country of origin* suggests that economics – especially in its theoretical mode – could have developed quite differently in the USSR too, had it not been banned and repressed on such a scale.

Building on several prominent cases and drawing on some previously unexplored material, this article examines the effect that the free exchange of ideas – or lack thereof – had on the productive coexistence of mathematics and economics. Apart from discussing missed opportunities, a more nuanced history also helps elucidate the role of the strong tradition in mathematics in the economists' theoretical innovation. The Soviet history demonstrates how this interdisciplinary interaction was shaped by the presence (or absence) of a diverse intellectual community as well as by the repressive political environment responsible for self-censorship and international isolation. Today, when academics, including economists, must confront politically influential ideologies – not only in Russia and in China, but across the globe – the historical reconstruction and contextualization of the Soviet experience may be valuable.

Although the primary focus of this paper is on Soviet mathematics, a brief description of the postwar Soviet *mathematical economics* is in order. By the end of the 1950s and early 1960s, the field had been institutionalized. This organizational process was driven by the mathematician Leonid Kantorovich and the statistician Vassily Nemchinov. With colleagues across various disciplines, Kantorovich and Nemchinov launched what was later–perhaps exaggeratedly—called the "mathematical revolution" in Soviet economics (Zauberman 1975a).

From the organizational point of view, the project was a success. After many battles with the Soviet Marxist orthodoxy, the new community of scholars succeeded in getting mathematics recognized as a legitimate tool in economic analysis. Consequently, multiple research centers were established across the country and university programs in "economic cybernetics" (an eclectic discipline that, in the Soviet context, could include a plethora of quantitative approaches to economic research and planning) became common. These programs produced several generations of scholars who were well-trained in mathematics and enthusiastic about its applications in increasing the efficiency of planning and industrial production.

---

[2] The similarity between Western and Soviet mathematics is clear from the multiple *parallel developments* more or less efficiently communicated through the Iron Curtain, and from the immediate post-Soviet history, when 'the […] diaspora of Soviet mathematicians […] led to both a mathematical labor supply shock and a mathematical knowledge shock in many countries' (Borjas and Doran 2013: 1157).

However, Soviet mathematical economics never became a part of global economic discourse. Isolated from the rest of the world, it remained exposed to ideological pressure inside the country (this made, for example, a detached study of the market economies hardly possible). Mathematical economists faced indifference, incomprehension, and outright hostility from the mainstream Soviet economics profession. They also lacked a well-defined theoretical framework in research and training.[3] This limited their ability to contribute directly to the development of economics in the way they could without these limitations. Their research often took a more abstract direction, using mathematical formalism to help hide particular policy messages such as the need to abandon inefficient overcentralized planning in favor of incentive-based, quasi-market mechanisms.

Another obvious limitation was that the Soviet economy was not ostensibly a market economy, the primary object of most economic analysis. When combined with the constraints on free academic conversation and the unavailability of data due to secrecy and censorship, this severely limited any attempts to study the workings of the Soviet planning system, too. As a result, the analysis of this system, initiated for example by the Hungarian economist Janos Kornai (1979; 1980), was unable to emerge as a standard topic of study in the Soviet intellectual landscape.[4]

The paper is structured as follows. Section 2 begins with a story of a conference in Moscow in 1935 that brought together many illustrious mathematicians. The narrative returns regularly to this event to explore its underappreciated relevance for the history of economics. In fact, the conference gathered, for the first *and* last time, most of my story's protagonists as well as their direct associates, who would later play various roles in the mathematization of economics in the US and in the Soviet Union. The section subsequently examines Soviet work on static and dynamic optimization and its significance for Soviet and Western economics. The discussion focuses on Leonid Kantorovich and sets the record straight on several events in his career that are missing or insufficiently elaborated in the existing narratives. Reassessing Kantorovich is necessary both in view of the new historical evidence and of the subsequent developments in economics. Tracing the work of another Soviet mathematician, Lev Pontryagin, and his colleagues will make a contrasting case of an influential contribution that was very quickly adopted in Western theories of economic dynamics, but not in Soviet mathematical economics.

Section 3 focuses on Soviet game theory. It will dwell on its general context and will illustrate the difficulties Soviet game theorists faced by drawing on the biography of Olga Bondareva, a talented scholar whose promising career was curtailed by the intervention of the Communist party apparatus and who could not become part of international game-theoretic community. This happened at precisely the moment when game theory was becoming a primary language of economic modeling – thus at a time when economic theory might have profited from the involvement of Soviet game theorists. Section 4 is devoted to probability theory and related work on the theory of stochastic processes in the Soviet Union. Besides discussing the ways this work influenced economics, it focuses on Andrei Kolmogorov—a classic in the field—and on his role in the mathematization of Soviet economics. Finally, it considers the various mathematicians coming from Kolmogorov's school. By tracing the import of their mathematical results, as well as the ways their engagement went beyond purely academic work, the section demonstrates the amazing inefficiency of interdisciplinary communication in contexts where doing economics was politically dubious and did not enjoy widespread societal acceptance. Section 5 concludes by reflecting on the lessons to be learned from this historical reassessment.

---

[3] Soviet mathematical economists, for the most part, were (self-) educated by the translations of Gale (1960/1963) and Karlin (1959/1964), but not by Samuelson (1947), Koopmans (1957), Dorfman et al. (1958), Debreu (1959) or even Arrow and Hahn (1971). While both the economic and the mathematical content of these two sets of books were very similar, the material that went beyond formalization and technical treatment of models was not allowed in the USSR and never got translated. The only exception was Samuelson's 1948 undergraduate textbook in economics. This edition, however, was plagued by all sorts of excisions, ideologically motivated editorial comments, and mistranslations (see Gershenkron 1978).

[4] However, there were notable exceptions to these broader tendencies, which should make us specifically appreciate the achievements of these scholars. (I am grateful to a referee for pressing me specifically on this point.) For the discussion of actors, institutions, and fields within Soviet mathematical economics, with some examples of these exceptions, see the Appendix.

## 2. New Mathematics for Economics on Both Sides of the Iron Curtain

### 2.1 Moscow, 1935: One-Time Encounters

In 1935, two young Soviet mathematicians, Pavel Alexandrov and Andrei Kolmogorov, brought together a remarkable group of scholars for the first international topological conference in Moscow. Topology, a new field of mathematics, dealing with a rigorous analysis of continuity, was on the rise, and Alexandrov was one of its active promoters (Alexandroff and Hopf 1935).[5] In hindsight, one can say that this gathering was even more important because it brought together American, European, and Soviet topologists just before the decades when direct personal contact was cut off for Soviet mathematicians.

Although political campaigns directed against mathematicians had begun in the 1920s, the situation deteriorated after 1935. Compared to the 1932 congress, no Soviet scholars participated in the 1936 International Congress of Mathematicians. While this had more to do rather with boycotting the Germans than forbidding to travel to 'capitalist countries,' the isolation only grew as the 1930s progressed. Fewer Soviet mathematicians published in Western journals (Alexandrov 1996) and international travel became increasingly rare. Moreover, in 1936 one of Moscow's senior prominent mathematicians, Nikolai Luzin, fell victim to a political campaign. He was criticized in the major Soviet newspaper, and forced to undergo the humiliating procedure of collective "criticism," in which his former students, including Alexandrov and Kolmogorov, participated. Ultimately Luzin was publicly declared an "enemy of the Soviet science" and lost his academic positions (Demidov and Levshin 2016).

Nonetheless, the 1935 topology conference had contributed to the increased visibility of Soviet mathematics, especially in the US, and led to formal and informal intellectual exchanges; for example, the work of some American mathematicians was published in the Soviet journals and vice versa. The parallelism of discoveries—something that reoccurs frequently in the narrative to follow—was also a hallmark of the conference. Most famously, it was quickly recognized during the discussions that both Kolmogorov and James Alexander had conceived, simultaneously and independently, the concept of cohomology, an important instrument in the study of topological spaces (Whitney 1988).

---

[5] Alexandrov's name first existed in German version (Alexandroff). The other protagonists' names are also sometimes spelled differently: Khinchin's name was initially spelled in French manner (hence the references to Khintchine in 4.1 below). Similarly, Nemchinov, Kantorovich, and Novozhilov could also be spelled differently, thus, Nemčinov, Kantorovič, and Novožilov (see Conclusion).

The most intriguing feature of the 1935 conference for the history of *economics* though was that it inadvertently gathered many of the scholars who would later play central roles in the mathematization of the discipline. The Soviet side included Kantorovich, Kolmogorov, and Pontryagin,[6] while among the American participants were John von Neumann and Albert W. Tucker. They were brought to the USSR with other US scholars by Moscow-born and Russian-speaking mathematician Solomon Lefschetz, Tucker's PhD supervisor. In the postwar years, Lefschetz would actively promote Pontryagin's work in the US mathematical community.

In the mid-1930s, however, the importance of advanced mathematical techniques for economic theory was not yet widely recognized. Among the Moscow conference attendees—and probably among anyone in the world—only von Neumann at that time could foresee that topological arguments (fixed-point theorem) would provide an entirely new formal framework for economics.[7] And of course, no one, Pontryagin included, knew that optimization principles could be applied to the control problems in dynamical systems; nor did anyone anticipate the effects this would have on modeling economic dynamics.

A decade later, the conversation was quite different. The new set of interconnected tools—known as linear and non-linear programming, "activity analysis," and game theory—anticipated by von Neumann's work (1928; 1937/1945), were be developed in the 1940s by George Dantzig, Tjalling Koopmans, but also, crucially, by Tucker and his colleagues (Harold Kuhn) and students (David Gale, John Nash).[8] By the end of the 1940s, one of Lefschetz's students, Richard Bellman, with support from von Neumann, started to develop dynamic programming that would prove a breakthrough in the field of optimal control.

The crux of the matter, however, was that unlike in 1935, such a conversation between Soviet and Western mathematicians *could not take place*. Soviet mathematics became increasingly isolated. Simultaneously, in Soviet academic economics any sort of quantitative reasoning was either absent or marginalized. The contrast between the ways mathematics and economics co-existed in the Soviet and Western worlds is striking. The cases discussed below demonstrate this in more detail.

## 2.2 Leonid Kantorovich and his (Partly Unknown) Legacies in Economics

Very early on, child prodigy Leonid Kantorovich (1912-1986) became one of the most promising Soviet scholars in pure mathematics. The mathematical theory that Kantorovich was supposed to present at the 1935 conference, and which was further developed by other attendees, explored a new class of partially ordered spaces introduced by Frigyes Riesz in 1928.[9]

---

[6] In fact, almost all the Soviet participants, including Israel Gelfand, Nikolai Luzin and Andrei Tikhonov, contributed, in various ways, to the mathematization of Soviet economics (Boldyrev 2024).
[7] Von Neumann first presented his general equilibrium model in 1932, completed the model in 1935, and published it in 1937. He cites Alexandroff and Hopf (1935) in his paper (von Neumann 1937/1945). At that time, the book was a standard reference in topology.
[8] See the historical accounts in Weintraub 1992; Myerson 1999; Giocoli 2003; Düppe and Weintraub 2014.

Surprisingly, even this, rather abstract, object turned out to be useful in economic modeling. When general equilibrium theory was extended to infinite dimensional commodity spaces (Bewley 1972; Aliprantis and Brown 1983), Riesz spaces (vector spaces with an additional ordering structure, also called vector lattices) became a helpful formalism.[10] Overall, this analysis is useful whenever continuous phenomena need to be modeled. These include, for example, intertemporal resource allocation, when consumption streams over infinite time horizon are considered; or equilibrium under uncertainty, when consumption is contingent upon the future states of the world, as in theoretical finance. The power of this approach became particularly clear once Harrison and Kreps (1979) and Kreps (1981) provided a new formalization of Black-Scholes-Merton theory.

But in the 1930s, Kantorovich was still only a mathematician. In 1938, he was approached by a group of engineers from the plywood trust, who wished to understand how to maximize their overall output, under certain technological constraints (Gardner 1990). He formalized this as a linear optimization problem, and in 1939 published a small book formulating the problems of this type in various planning contexts. The book also detailed the solutions to such problems and sketched an existence proof of the solutions. Subsequently, Kantorovich migrated to economics; he did not do any mathematical research after the 1950s.

Kantorovich's (1939/1960) work belonged to a broader project of restating the economic problem–the optimal allocation of scarce resources–in a more rigorous way. He approached it from the perspective of production planning. Importantly, however, his approach, like that of von Neumann, was informed by his interest in functional analysis. The fundamental techniques both Kantorovich and von Neumann used, based on the geometrical versions of the Hahn-Banach theorem, involved the so-called separating hyperplane arguments. The idea of a duality between technical/production space and the space of monetary units – or, as economists would say, between the "real" and the "nominal"–is explicit in von Neumann's (1937/1945) influential model of optimal growth as well as implicit in Kantorovich (1939/1960).[11] Tellingly, the genealogy of these contributions was different. Von Neumann came to his model after experimenting with various other modeling frameworks and exploring the economic equilibrium theories available to him at the time–namely, those of Walras, Wicksell, and Cassel (Carvajalino 2021). In contrast, Kantorovich arrived to his results by formalizing an applied planning problem.

---

[9] Kantorovich did not actually give a talk at the conference because he could not properly prepare it on time due to illness. The most complete presentation is in Kantorovich 1937. He sent the paper to von Neumann in 1936 and received a positive response, but did not publish it in *Annals of Mathematics*, apparently because von Neumann wanted to shorten the paper, see the details in Boldyrev 2024.

[10] This structure implies that every finite nonempty subset of a given vector space has the supremum – least upper bound – and infimum – greatest lower bound (Aliprantis and Burkinshaw 1978). Curiously, Kantorovich (1966) argued that the order structure of vector spaces could be used for the analysis of preferences in economics but did not go beyond this vague indication. Since the 1980s, lattice theory has played an increasingly important role in mathematical economics toolbox.

[11] Note that Kantorovich did not use a fixed-point theorem to demonstrate the existence of what he later called "objectively determined valuations," or shadow prices. Instead, he made recourse to the simpler geometric argument. The geometric version of the hyperplane separation theorem is due to Minkowski (1911). See also the historical account in Kjeldsen 2001. At the 1985 World Congress of the Econometric Society at MIT, Debreu (1986: 1261) said that Minkowski's theorem had to wait more than thirty years until applied to economics. Kantorovich, who was, curiously, among the audience of Debreu's lecture in Cambridge, MA, could have corrected him.

The economic meaning of this contribution was well captured by Hicks (1960: 707) who surveyed the new work in mathematical economics:

> It has been made apparent, not only that a price system is inherent in the problem of maximising production from given resources – but also that something like a price system is inherent in any problem of maximisation against restraints. The imputation of prices (or "scarcities") to the factors of production is nothing else but a measurement of the intensities of the restraints; such intensities are always implicit – the special property of a competitive system is that it brings them out and makes them visible. It is through its power of developing the intensities (in the photographic sense of developing), so that they are available for use as instruments in the process of maximisation, that the competitive system does its job.

The discovery of the role prices play in rational economic activity, emphasized informally by Mises (1920), was something Kantorovich himself considered to be his major contribution to economics and which, arguably, brought him the Nobel Memorial prize he shared with Koopmans in 1975.[12] And it was precisely this contribution that made his subsequent career as an economist in the USSR so difficult (Boldyrev and Düppe 2020). The actual application of linear optimization would require the deep reform of both Soviet economic science and the planning practices– something hardly possible at the time.

At first, Kantorovich was very enthusiastic. As early as 1942, he had prepared a book-length exposition of linear theory and its economic applications. However, his multiple attempts to gain the attention of economists and policymakers faced a fierce ideological resistance. Mathematical methods, including Kantorovich's methods of pricing relative scarcities, were deemed ideologically suspect. Only after 1956, supported by his colleagues in economics and mathematics (see also Section 4.1 below), did Kantorovich build up a new community (on the vicissitudes of this institutional process and on the resistance he had to face see Boldyrev and Düppe 2020).

During World War II, Kantorovich published a purely formal paper, "On the translocation of masses" (Kantorovich 1942), in which he provided a treatment of the optimal transportation problem. The problem, first formulated by Gaspard Monge in the eighteenth century, consists in finding the optimal transport plan – for moving the resources across economic units – that would minimize the cost of transportation. Kantorovich formalized this problem in terms of probability measures in a metric space. In some sense, Kantorovich's solution was an infinite dimensional version of the duality theorem. Astonishingly, the ideas from this paper (see its history in Vershik 2013) and from the subsequent research done by Kantorovich, his student Gennady Rubinstein (Kantorovich 1948; Kantorovich and Rubinstein 1958), and others, have only recently been taken up in advanced formal economic theory. Optimal transport techniques have found applications in multiple fields ranging from asymmetric information to labor economics to theoretical econometrics (Ekeland 2010; Galichon 2016; 2021).

It was this paper, inspired by the actual planning problems, that attracted the attention of Koopmans (who also worked in the field) and convinced him of Kantorovich's priority in the development of linear programming. The paper had been published in English in 1958, even before Koopmans organized the English translation of the first, economic paper (Kantorovich 1939/1960).

---

[12] Arrow nominated Kantorovich for the Nobel prize in 1972 (see Box 56, Kenneth J. Arrow Papers, David M. Rubenstein Rare Book & Manuscript Library, Duke University).

By that time, operations research in the US was already in full swing. Kantorovich's 1942 book would only appear in Russian in 1959 and in English in 1965 (Kantorovich 1959/1965). But while the American economics profession quickly absorbed the new ideas (Dorfman et al. 1958), there was no comparable Soviet economics profession to take it up on the other side of the Iron Curtain.

Ironically, Kantorovich received credit not as a mathematician, but as an economist. This was largely because linear theory was developed nearly simultaneously in the West and in the East and the discussion of priority was unsettled (Dorfman 1984; Gass 1989; Bockman and Bernstein 2008). Nevertheless, the economists Koopmans and Dorfman believed in the priority of Kantorovich precisely because of the economic interpretation he gave to what he called "resolving multipliers" (that is, shadow prices emerging in any linear optimization problem). However, operations researchers considered linear programming a part of applied mathematics ("the body of theory and methods to solve linear programs with computational efficiency" (Schwartz 1989: 146)). They argued that Kantorovich (1939/1960) did not prove a duality theorem and did not provide an efficient method to solve the problem – the two elements necessary to constitute a true discovery (Gardner 1990).

Wherever the priority lies, the actual global *influence* of Kantorovich's linear programming work in economics has been limited, at least formally. In the West, he was long the least cited of all economics Nobelists (Bjork et al., 2014). This was not true, however, in the USSR and East European countries. Kantorovich's work and organizational efforts were decisive for institutionalizing Soviet mathematical economics and for developing Soviet optimization theory (on the latter see the overview in Polyak 2002).

The parallels between Soviet and Western work help illuminate differences, which are particularly salient in the interpretation of the same mathematical results. While for Koopmans, the invention of linear programming meant that the methods for managing scarcities are universal for all economic systems, Kantorovich repeatedly insisted that socialist planning was the form of economic organization uniquely suited for applying optimization methods at a large scale. In fact, while remaining loyal to the regime throughout his career and masterfully drawing on the ideological and bureaucratic resources of the Soviet academic system, Kantorovich used this argument strategically, to convince the policymakers that mathematics matters.

## 2.3 Lev Pontryagin: Cold War Mathematics and Economic Dynamics

Pontryagin's name is also quite familiar to many economists. But, unlike Kantorovich, he never became interested in economics. At the Moscow conference in 1935, he presented his work on topological groups which was closely connected to von Neumann's research at that time. But by the 1950s, following the suggestions of the academic bureaucracy (Pontryagin 1998) or perhaps discouraged by the superior results of the French topological school (Boltyansky 1999), Pontryagin decided to abandon topology altogether and to turn to applied mathematics.

While Kantorovich's work was motivated by planning problems, the formulation of what is now called Pontryagin's "maximum principle" was developed for military research. The Soviet military was curious whether the processes of aircraft pursuit and evasion could be described mathematically. In attempting to solve this problem, Pontryagin and his group, which included Vladimir Boltyansky, Revaz

Gamkrelidze, and Evgeny Mishchenko, discovered the need to re-invent another field of mathematics: the calculus of variations (Gamkrelidze 1999).

As with the development of linear programming in the 1930s-1940s, many advances occurred in parallel. In the field of dynamic optimization, the RAND corporation scholars Magnus Hestenes, Rufus Isaacs, and Richard Bellman were pioneers. Pontryagin and his group began working on optimal control a bit later in the 1950s (Pesch and Plail 2009). The similarity in motivations between East and West in expanding the calculus of variations was also remarkable: both groups of scholars were inspired by engineering and military applications. In 1955, Pontryagin's group managed to prove the "maximum principle" in its general form. The ultimate set of results was published in 1961. It was very quickly translated into English (Pontryagin et al. 1962) and almost immediately adopted in economics.[13] Although other Western economists (Ramsey 1928; Samuelson and Solow 1956) had previously used the calculus of variations in growth theory, it was the formalism suggested by Pontryagin et al. that became most popular in economic dynamics.

Similar to linear programming, the calculus of variations helps illustrate the diverging interpretations given to mathematical results. In his Nobel lecture, Paul Samuelson (1972) argued that to get a neater, "a better, a more economical, description of economic behavior" (p. 251) one needed the idea of maximization. For Samuelson, the description of Newton's falling apple in the way that "its position as a function of time follows that arc which minimizes" a certain integral (in other words, describing it in the framework of the calculus of variations) is as useful as describing the behavior of economic agents in terms of optimization under constraints. And here, the comparison with Kantorovich is instructive. In a talk given in 1960, he argued that economic laws in socialist societies get actualized and implemented by planning decisions. Curiously, Kantorovich, too, invoked the same mathematics, but with a different, more "interventionist" overtone:

> There are two classical problems [in the calculus of variations]. One is about a catenary, and the other is about the curve of the fastest descent. The chain sags along the catenary line regardless of whether the person who holds its ends knows calculus of variations or not. But whether the descent will be along the desired curve – the cycloid – depends on the extent to which the designer knows calculus of variations (Kantorovich et al. 2004: p. 119).[14]

After the English publication of Pontryagin et al. (1962), young Stanford economists Hirofumi Uzawa, David Cass, and Karl Shell, informed by their colleagues Samuel Karlin and Kenneth Arrow, applied the maximum principle to optimal growth theory. Pontryagin et al.'s work on the maximum principle gave the growth theorists straightforward methods to put "the transformation of the equations of motion into dual form that has an important economic interpretation" (Magill 1970: 13). The central instrument in Pontryagin et al. – namely, the so-called Hamiltonian function— was interpreted as the maximized value of an aggregated economic indicator (Shell 1987). The popularity of the technique, which extended the principles of constrained optimization to a dynamic case, went so far that the maximum principle was heralded

---

[13] The reception was quick. Burmeister (2009: 38) recollects that Solow, in his lectures, applied Pontryagin's "minimum principle" to the solution of Ramsey's problem as early as December 1962. Uzawa (1965) first submitted his paper in 1963, in a year after Pontryagin et al. (1962) first appeared in English.

[14] Samuelson was reluctant to think of the economic equilibrium analysis in terms of planning (Hands 1994).

"to be the culmination of a logical sequence originating in the maximum principle of Adam Smith" (Burmeister and Dobell 1970: 404, cited in Wulwick 1995: 421).[15]

Importantly for my argument here, there was no way such an economic interpretation – and the surrounding enthusiasm – could have emerged in the Soviet academy. Indeed, Pontryagin, a loyal communist, would have been the first to condemn it strongly and emphatically. This helps account for the fact that the application of the maximum principle in Soviet mathematical economics was smaller in scope, happened later than in the West, and remained quite piecemeal (Zauberman 1976). This did not change even after Cass (1965) and several other subsequent optimal growth models were framed as models of a centralized planned economy featuring a "social planner" or a "central planning board."

This literature was read in the Soviet Union. In 1971, the work of American economists on optimal control and economic growth was being taught by Boris Mityagin (1972) at Moscow State University in the mechanics and mathematics department (*mekhmat*). While it found its natural audience there, the delay in reception was noticeable.

## 3. Episodes in the History of Soviet Game Theory

### 3.1 The General Context: The Games Soviets Played

We do not know whether von Neumann ever discussed game theory with his fellow mathematicians in Moscow in 1935. We do know, however, that a year later one of the conference participants, Eduard Čech, served as a mediator between von Neumann and Oscar Morgenstern, effectively enabling their now famous collaboration (Morgenstern 1976, 806f.). Another participant, Israel Gelfand, later helped initiate the research on games of finite automata. And in the 1950s yet another participant, Andrei Markov Jr., was the supervisor of Nikolai Vorob'ev, the founder of Soviet game theory.

For decades, Vorob'ev was an active researcher and popularizer of this field in the USSR. His first short expository paper on game theory appeared as early as 1955. Most likely, his interest was sparked by his engagement with probability and statistical inference, which at that time was formalized by Wald (1949) as a game with Nature. In fact, the first ever book-length treatment of game theory in the USSR was the translation of Blackwell and Girshik's *Theory of Games and Statistical Decisions* (1954/1958), which develops Wald's approach. Other publications quickly followed. While early game theory in the United States was primarily applied in military contexts and in international relations (Mirowski 2002; Erickson 2015), the Soviet military, too, were paying attention. Indeed, an anthology on the military applications of game-theoretic models (Ashkenazi 1961) appeared among the first in this stream of literature.[16]

Vorob'ev's efforts were tireless. He supervised PhD students, wrote review papers and a textbook in game theory (Vorob'ev 1977), and organized All-Union game-theoretic conferences (1968, 1971, 1974). These activities helped to establish a sizeable research community, primarily in Leningrad (Hagemann et al. 2016). The members of this community included, in different time periods, Olga Bondareva, Victor Domanski, Victoria Kreps, Tatiana Kulakovskaya, Natalia Naumova, Iosif Romanovskii, Arkady Sobolev, and Elena Yanovskaya, among others. While the initial focus was on developing algorithms for computing equilibria in non-cooperative games (Vorob'ev 1958; Romanovskii 1962; Raghavan 2002; Stengel 2002), later, it shifted to modeling cooperative interactions (many of the results are summarized in Maschler 1992). By 1966, when Morgenstern visited Moscow for the International Congress of Mathematicians, Soviet game theory was already firmly established. In 1971, when the *International Journal of Game Theory* was launched, Vorob'ev sat on its editorial board.

---

[15] By 1970 the technique became standard in economics as seen in Arrow and Kurz (1970).
[16] Perhaps the most prominent Moscow scholars in operations research at that time, David Yudin and Evgenii Gol'shtein (also credited as Golstein), were working at a classified military research institute (NII-5). Yudin and Vorob'ev edited the Russian translations of major game-theoretic literature. See also the overview of early Soviet ideological debates around the applications of game theory in Robinson 1970.

Outside Leningrad, in the 1970s and 1980s, a lot of game-theoretic research was done in Moscow. Yuri Germeyer and several of his colleagues were working on hierarchical games and trying to formalize information asymmetries (Germeyer had connections to military research organizations as well). Furthermore, Michael Tsetlin and Gelfand were doing pioneering research on the games played by finite automata.[17] Somewhat later, Vladimir Danilov (1992) at CEMI, the foremost institution in Soviet mathematical economics (see Appendix) was working on Nash implementation. This work contributed to the larger mechanism design literature investigating the properties of social choice rules, under which the set of the optimal outcomes would coincide with the set of Nash equilibria. Game theorists were also active in other Soviet regions. In Novosibirsk, Valeri Vasil'ev was applying the theory of vector lattices to cooperative games and a group led by Alexander Granberg was building game-theoretic models in spatial economics. In Vilnius, Vorob'ev's first student Eduardas Vilkas and his associates were also working on cooperative games.

It was natural at that time to consider game theory as simply part of mathematics. Most work in the USSR was published in mathematical journals and was rarely motivated by economic problems. At the same time, real planning involved conflicts and strategic behavior. The popular decomposition algorithms for solving linear problems, for example, featured a game between the planning center and the planned sectors of the economy (Kornai and Liptak 1965). But while the official economic ideology had some room for "optimal planning," it would never endorse the modeling of conflicts between self-interested agents. Hence, the attempts to apply game theory in Soviet planning were too few and too inconclusive to claim that it could inform policymaking (Zauberman 1975b). The same was true for any other aspect of social life. But even granting the abstractness and political neutrality of their discipline, the general context of Soviet science was stifling for game theorists, as will be illustrated in the next section.

### 3.2    Olga Bondareva: The Missing Coalitions

The story of Olga Bondareva (1937-1991), a game theorist from Leningrad, is instructive for understanding both the character of Soviet applied mathematics and the constraints imposed on its development by the ideological environment.

---

[17] Gelfand (1913-2009) was one of the most influential Soviet mathematicians. In economics, the idea of modeling games with agents characterized by simple behaviors was then taken up much later. It was developed, in various contexts and following various theoretical traditions, by Rubinstein (1986) and Roth and Erev (1995). At stake was the integration of bounded rationality and learning into game theory.

Bondareva's contribution to game theory was both pathbreaking and influential. Working under the supervision of Vorob'ev, Bondareva (1962; 1963/1968) provided the first characterization of cooperative games with a non-empty core, a concept analogous to equilibrium/optimality. She managed to attach numbers, also called balancing weights, to the partition of the set of players (a collection of this set's non-empty subsets, such that each player belongs to only one of these subsets). In a world of balanced collections of players, a simple linear inequality helps check whether the core is non-empty. This meant that checking the existence of the non-trivial core became easier, and this condition was essentially connected with linear optimization. In fact, Bondareva, as well as Shapley (1967) who arrived at the same result independently, refer to the duality arguments in linear programming. This result was very quickly noticed in the game-theoretic community and became a classical theorem in cooperative game theory, with many economic applications (Wooders 2008).

Bondareva's promising career at the mathematics and mechanics department of Leningrad University (*matmekh*) changed direction quite dramatically. In 1971, a student who attended Bondareva's game theory seminar applied to emigrate to Israel, which, at that time, was the only legal way to leave the country. Following routine procedure, he had to be expelled from *Komsomol*, the Communist youth organization. At the Komsomol meeting that considered his case, Bondareva's quite conformist argument was not to defend the student's decision, but to let him stay in Komsomol "out of pity," because leaving one's homeland was already a sort of a suicide. But then she became more critical; her indignation was spurred by the general antisemitic sentiment of the whole discussion. In the end, the vote to expel the student was not unanimous and the authorities blamed Bondareva's intervention.[18] Bondareva was required to leave the university and, only after some time, was she able to find a job in the Leningrad University economics department.

For Bondareva the academic, this was difficult. Economics was never in high regard among her mathematician peers. The entire economics department was deemed "ideological," and her ambitions were elsewhere. Worse still, due in part to her strained relations with Vorob'ev, she was banned from the Leningrad game theory community and could continue her work only by informally meeting with former students and colleagues.

To fully understand this episode, one should account for the gender aspects in Bondareva's non-obedience story. Her husband recollects that a male colleague, around the same time, did something worse than her: he wrote an anonymous letter to the Party's Central Committee condemning the 1968 Soviet invasion in Czechoslovakia (Gordon 1992). Nonetheless, his punishment was less severe. He was suspended but not fired because the university bureaucracy held that his letter was just a delusion rather than a sign of a systemic disloyalty.

Another important dimension in this story was the sheer number of women game theorists in Leningrad, many of whom were either trained by Bondareva or were encouraged by her work and career. At that time, this was not a typical career pattern. While not formally disadvantaged, women (some prominent counterexamples, like Olga Ladyzhenskaya, notwithstanding) mostly chose not to engage in mathematics, pure or applied.

---

[18] The character of late Soviet public life (Yurchak 2005) seems to have played a role here. Party bureaucrats could have supported Bondareva and even share her beliefs, but the public expression of outrage and non-obedience on her part, especially if it made other colleagues rethink their views, was unacceptable.

During some of her interim research appointments, Bondareva was actually confronted with economic material (this involved contract work on solving linear programming problems). However, she would never embark on broader economic applications, and probably never wished this, likely because game theory, in her mind, was a mathematical discipline, not a tool for economic analysis. Still, soon after she began teaching in the economics department, she published a short book on game theory and its applications to economics (Bondareva 1974). The book mostly focused on cooperative games and considered several real-life situations that could be elucidated via game-theoretic solutions: a market with three agents, voting schemes, coalitions of workers deciding on the payment for a job, inventory policy, and optimal routes. However, in Bondareva's narrative these problems serve mostly as illustrations. "Unfortunately, she writes, "we are not aware yet of any serious enough economic problems that were solved with this theory (of cooperative and non-cooperative $n$-person games – I.B.), but this situation might change in the nearest future" (Bondareva 1974: 37).

The timing is crucial here. Bondareva's conjecture proved right. The 1970s was a period of intense work in game theory, resulting in a profound transformation of economics (documented for example in Myerson 1999). Even previously, the concept of Shapley value, another key idea from cooperative game theory (simply put, a reasonable and axiomatically justified way to allocate individual payoffs in these games), had been applied in political science (Shapley and Shubik 1954). Most members of the Leningrad group were working with this concept, but given the absence of a political science discipline and of a democratic political system, similar applications of game theory were inconceivable in the Soviet Union (see, however, Boldyrev 2020 on an exception).

Bondareva's academic isolation within the country made the effect of the Iron Curtain even worse. In the years following the international recognition of her result, she was not permitted to travel abroad nor connect with international colleagues, despite multiple invitations. Further, Bondareva was never able to reach the top of the Soviet academic hierarchy. She received a professorship only near the end of her life (Kulakovskaya and Naumova 1992). These last years (1988-1991), before she suddenly died in a car accident, were marked by fascinating collaborations, and the realization, on Bondareva's part, that for decades there was a scholarly community, a coalition of sorts, for which her work was important, even if she could not be a full-fledged interlocutor.[19] As Roger Myerson stated in his obituary, "[w]ithout direct personal contact, our knowledge of each others' work was necessarily limited, to the detriment of scientific progress in both East and West" (1992: 324).

What applied to Bondareva, was true for her colleagues as well. They considered game theory to be a part of mathematics, and they mostly avoided any economic applications. Of course, the obvious and popular applications of game-theoretic reasoning at that time (industrial organization, auctions, matching or voting) were not feasible topics in the Soviet Union. And even in the West, where this adverse ideological pressure was not a factor, game-theoretic modeling in *economics* only became widespread at the beginning of the 1980s, after the zenith of Soviet mathematical economics. Thus, despite the relevance of their mathematical research,

---

[19] In 1989, when Ehud Kalai founded the journal *Games and Economic Behavior,* almost all the major game theorists of the time were asked to serve on its editorial board. Bondareva was the only Soviet and the only woman in this group of 47 scholars. At the time of writing, eight of these individuals have been awarded the Nobel Memorial prize in economics.

Soviet game theorists were not able to participate in the broader movement that eventually made game theory the formal language of economics.

## 4. Soviet School of Probability Theory: Contributions in Context

### 4.1 Kolmogorov and Economics: A Story of Inspiration and Patronage

Alexandrov co-organized the 1935 Moscow conference with his partner, Kolmogorov. At that time, Kolmogorov was already known for his axiomatic restatement of probability theory (1933), which became central for the whole field and defined the language used to conceptualize probability in contemporary mathematics.

Kolmogorov was not alone in revolutionizing probability theory. Generally, Soviet scholars built on important legacies from pre-revolutionary Russia, both in mathematics (in the work of Pafnuty Chebyshev and Andrei Markov) and in statistics (where Alexander A. Chuprov was a key contributor). These legacies continued to influence the research and teaching cultures in the key academic centers of Moscow and Leningrad, as well as beyond. Perhaps the most well-known of Kolmogorov's immediate interlocutors was Aleksandr Khinchin (1894-1959). While Kolmogorov's work is a separate chapter in the history of science, for the story here two aspects are particularly important.

*First*, Kolmogorov's own contributions in various fields related to probability and stochastic processes have inspired many economic applications. For example, in the 1930s, Kolmogorov's (1933) and Khinchin's (1933, 1934) formalizations of stationary stochastic processes helped the Swedish statistician Herman Wold (1938) create a full-fledged theoretical framework for the analysis of stationary time series (see Mills 2011). Their work on queuing theory (Kolmogorov 1931a; Khinchin 1932; Khintchine 1960), in which Khinchin provided some basic formalisms, was clearly motivated by economic considerations and applications.[20] While integrating the concept of entropy into probability theory, Khinchin (1953, 1957) provided axioms for Shannon's entropy index that could be used as a measure of heterogeneity in populations. This work thus paved the way for applying probability theory to the analysis of social and economic diversity and inequality that is becoming increasingly important in economics.[21]

Kolmogorov's results in probability and stochastic processes have found numerous other economic applications, of which only a few can be discussed here. For example, the Daniell-Kolmogorov extension theorem, which constructs a probability space and a stochastic process from the collection of finite-dimensional distributions, has become an important instrument in the theoretical analysis of time series (Brockwell and Davis 1990). A statistical goodness-of-fit tests invented in the 1930s by Kolmogorov and Nikolai Smirnov was introduced to economics via statistics in the 1950s (Massey 1951).

With the advent of microfoundations and rational expectations, the task of estimating dynamic macroeconomic models came to involve agents facing stochastic processes. The parameters of these processes had to be disentangled from the parameters of the agents' objective functions, which "would enable the econometrician to predict how agents' decision rules would change across alterations in their stochastic environment", or their constraints (Hansen and Sargent 1980: 7-8). It was in deriving formulas for decision rules that Wiener-Kolmogorov prediction theory came to be used by these authors.[22] A related problem, that of filtering, or, simply put, purifying random processes from noise, was also addressed by Wiener's and Kolmogorov's work. Wiener-Kolmogorov filtering is now a widely applied technique in time series analysis.

---

[20] In the Soviet literature, this applied field of probabilistic modeling had a more politically correct name, the "theory of mass service." After WWII, the Soviet scholars who suffered the very real problems of shortages and queues, did not adopt the new term.

[21] Khinchin's axioms were later simplified by a Leningrad mathematician, Dmitry Faddeev (1956). I am grateful to a referee for this reference to Khinchin's and Faddeev's work.

Finally, Kolmogorov's (1931b) characterization of continuous stochastic processes via what came to be called Kolmogorov equations has also become part of the mathematical "infrastructure" informing economic analysis.[23] When, for example, macroeconomic shocks need to be modeled as diffusions similar to Brownian motion, this early theory is still indispensable (Dixit and Pindyck 1994; Stokey 2009).

But there was a *second* aspect to Kolmogorov's role in economics that is less known and rarely appreciated. Kolmogorov had deep interest in all sorts of applied mathematics. His first academic paper was in what is now called cliometrics, or quantitative economic history. In the second half of the 1950s, when cybernetics in the USSR ceased to be "bourgeois pseudo-science," as it was commonly chastised, and became a legitimate field of study, Kolmogorov naturally became part of this movement. There was strong enthusiasm for applying mathematical ideas to broader contexts and Kolmogorov saw the mathematization of economics as one of the key applications.

Kolmogorov's interest in economics could also have been inspired by Eugene Slutsky (1880-1948), a former member of the Conjuncture Institute and one of the most significant Russian economists of his time. Slutsky, in fact, began his academic career before the Bolsheviks took power and worked in both mathematical economics and statistics.[24] Given the worsening political climate, by the end of the 1920s Slutsky had abandoned economics and turned to statistical modeling in meteorology and to the theory of stochastic processes. Myerson's (1992) words of regret cited in the previous section echo those by R. G. D. Allen (1950: 213f.), written at the beginning of Cold War:

> It is unfortunate that, for so long before his death, Slutsky was almost inaccessible to economists and statisticians outside Russia. He opened up new areas but left them to be explored by others, and the exploration even now is far from complete. His assistance, or at least personal contacts with him, would have been invaluable.

In an obituary, Kolmogorov (1948: 147), whose own work in probability theory owed much to Slutsky (Shafer and Vovk 2006), wrote that Slutsky's research in natural science and economics "is interesting primarily from the methodological point of view as indications of ways to expand the power and flexibility of mathematical methods." But Kolmogorov surely knew that Slutsky's engagement with economics, including his famous contributions to demand theory (Slutsky 1915) and business cycle research (Slutzky 1937), was anything but accidental or piecemeal. The same can be said about Kolmogorov himself.[25]

---

[22] The theory comes from Wiener (1949) and Kolmogorov (1941 a, b). They developed it essentially in parallel (see historical comments in Wiener 1949: v, 59). The ideas were communicated to Sargent and other macroeconomists via Whittle (1963). The book was so useful that Sargent wrote a preface to its 2nd edition in 1983. For a detailed account of the historical context, see Sent (1998).

[23] More precisely, these are differential equations helping describe the probability distributions of a state variable of a continuous-time Markov process.

[24] Sliutsky also taught economics. One of his students in Kyiv in 1915-1916 was Jacob Marschak, the future head of the Cowles Commission. Marschak later described Sliutsky as "one of the founders of mathematical economics" (cit. in Barnett 2011: 25).

Over the years, Kolmogorov's authority in Soviet science only grew. In 1939, he became a full member of Soviet Academy of Sciences, the highest position for a scholar in the USSR. In this capacity, he actively supported the mathematization of economics. Early on, he endorsed Kantorovich's applied study on optimal transportation; in 1957, he insisted on the need to publish Kantorovich's book manuscript, which he read closely. In 1960, when Kantorovich and Nemchinov organized a huge conference inaugurating the universal acceptance of "mathematical methods" as a legitimate tool in economic analysis, the intervention of Kolmogorov was decisive. Indeed, he interceded at the crucial point of contention, ideologically defending Kantorovich's approach to pricing the scarcities by providing a suitable interpretation of "resolving multipliers" (see the section 2.2.). In 1964, Kolmogorov criticized Kantorovich's critics who claimed Kantorovich should not be awarded the prestigious Lenin prize. Kolmogorov argued:

> Multipliers appear when solving problems for conditional extrema in a wide variety of fields. If we declare them, in advance, to be harmful indicators of "marginal utility," we will not only have to reject the work of L. V. Kantorovich, but also to declare the mathematical discipline called *calculus of variations* to be anti-Marxist. (Kolmogorov in Kantorovich et al., 2004: 335).

Another statement in this letter reads as an indictment of Soviet economic science, emphasizing the need

> to bring the price system into such a state, in which the economic calculation of the enterprise's monetary profitability would reflect a real approximation to satisfying the needs of society. As far as I know, no other general theoretical solution to the question of how to establish such a price system has been proposed except for theories equivalent to L. V. Kantorovich's theory of objectively determined valuations. (Kolmogorov in: Kantorovich et al., 2004: 335).

Finally, as a full member of the Academy of Sciences, Kolmogorov had the authority to communicate the papers to its Proceedings (*Doklady Akademii Nauk*), a prestigious journal and the central outlet for Soviet science in all disciplines. In the 1950s and 60s, most of the papers on linear programming, game theory, and related fields of mathematical economics were communicated by Kolmogorov, effectively bringing formal economic models into the legitimate academic conversation.

Unlike Kantorovich, however, Kolmogorov himself never became an economist, and did not seem to care about the applications of his ideas in Western economics. Nevertheless, while his support for rigorous economic analysis was mostly institutional and ideological, it was no less real.

---

[25] Slutsky defended his thesis *Theories of Marginal Utility* in 1910, and his last paper on the foundations of economics, a critique of Böhm-Bawerk's value theory, was published in German in 1927 (his 1937 *Econometrica* paper also dates back to 1927). Kolmogorov's 1957 letter to Kantorovich sharing comments regarding his book manuscript demonstrated that Kolmogorov knew of the concept of marginal utility. This was atypical for a Soviet mathematician at that time. The concept itself was not mentioned in Kantorovich's manuscript (in part due to political self-censorship) and would only be very slowly and cautiously adopted by Soviet economists. The fear of "dirty words" was pervasive. For example, Viktor Volkonskii (1967), a gifted mathematician and a student of Eugene Dynkin (see 4.2 and Appendix), who made one of the most sustained attempts to integrate the theory of optimal planning with general equilibrium analysis, would still call *u* a "preference function," rather than utility, would refer to it as a consumer valuation etc.

## 4.2 Soviet Probabilists and Economics: Imaginary and Real Encounters

In the 1960s, Kolmogorov's school of probability theory and stochastic processes became an important hub of new ideas during the "Golden Years" of Moscow mathematics (Zdravkovska and Duren 2007). The work on stochastic processes, including Soviet contributions, was summarized in book-length expositions that were quickly translated into English and became widely read in economics.[26] Exemplary contributions came from Kolmogorov's students Yuri Rozanov, Anatoly Skorokhod, and Roland Dobrushin.

Christopher Sims (1972) used definitions and results from Rozanov's (1967) book when developing the idea of Granger causality and attacking a difficult problem of the exogeneity of money in distributed lag regressions. Money and GNP were modeled as stochastic processes. In the 1950s, Skorokhod (1965) formulated techniques that allowed him to prove limit theorems for the broader classes of stochastic processes, including discontinuous ones.[27] These techniques were applied e.g., in the influential econometric analysis of unit roots in time series (Phillips and Perron 1988). The work of Dobrushin (1968) on generalized stochastic processes (random fields) informed the research on the systems of interacting agents and generally on the new probabilistic representations of the economy (Föllmer 1974; Durlauf 1993).

In these cases, like in the case of Kolmogorov, a series of new technical results turned out to be important for various fields of economics outside the Soviet Union, while the contributors themselves did not participate in these economic applications. Probably they were unaware of them: economic modeling was, for the most part, not on their research agenda.

An interesting example of a scholar from the same generation, whose work was also impactful in mathematics and related fields, but who turned to economics, is Igor Girsanov (1934-1967). Girsanov's (1960) theorem was widely applied in mathematical finance allowing modelers to change probability measures from the "observed" to the "risk-neutral" ones.[28] However, this result was only integrated into finance after Harrison and Kreps (1979) applied stochastic calculus (martingales theory) to asset pricing.

---

[26] See, for example, an introductory book by Yaglom 1962 or a more advanced treatment in Rozanov 1967 and Gikhman and Skorokhod 1969. Apart from the translations, an important role in the dissemination of this work was played by Billingsley's (1968) codification of the topic, which also covered the work of other Soviet probability theorists.

[27] Skorokhod's other supervisor was Kolmogorov's student Eugene Dynkin. Dynkin was dismissed from Moscow University in 1968 after signing a letter in support of the dissidents, and found a refuge in the CEMI, where he worked on economics together with a group of other notable mathematicians (see Appendix).

[28] As Braumann (2019: 6) helpfully explains, '[t]he idea in such applications [of the Girsanov theorem – I.B.] with risky financial assets is to change its drift to that of a riskless asset. This basically amounts to changing the risky asset average rate of return so that it becomes equal to the rate of return $r$ of a riskless asset. Girsanov theorem shows that you can do this by artificially replacing the true probabilities of the different market histories by new probabilities (not the true ones) given by a probability measure called the equivalent martingale measure. In that way, if you discount the risky asset by the discount rate $r$, it becomes a martingale (a concept akin to a fair game) with respect to the new probability measure. Martingales have nice properties, and you can compute easily things concerning the risky and derivative assets that interest you.'

Girsanov's short academic career is less well known. Inspired by the perceived potential of applying "mathematical methods" in socialist planning, he became involved in promoting the mathematization of economics in the USSR. Girsanov organized a group of lecturers who taught mathematics at the Moscow University's economics department and assisted Nemchinov in defining a mathematical curriculum for economics students (see Appendix). Furthermore, Girsanov began working on optimization theory, he gave several talks at the economic conferences, prepared a seminar (*spetskurs*) on optimization at the *mekhmat* of Moscow University in 1963 and 1964, and even taught linear programming at a high school.[29] He thus seemed to be very much willing to contribute to the new mathematical economics community.

But the community Girsanov was helping to build was not able to absorb and put to work his–and others'–most innovative mathematical ideas.[30] Macroeconomics could not be openly discussed, because economic policies, especially since the end of the 1960s, were not subject to critical public or academic deliberation. Inflation was officially non-existent. As far as econometrics was concerned, its development was clearly hampered by the unreliability and limitations of Soviet statistics and data gathering. The very problem of asset pricing, for which the Girsanov theorem and related mathematical ideas became so useful, could not be really discussed and indeed, seemed entirely irrelevant in a country where no financial markets operated until 1991. For example, it was not until the early 1990s that Yuri Kabanov, the young CEMI mathematician focusing on stochastic calculus, began to work on mathematical finance (see his recollection in Appendix).

## 5. Conclusion

What general lessons can we draw from this historical exploration? Apart from enriching our understanding of the past century's economics, and particularly of the ways it became mathematized, this narrative illustrates how interdisciplinary entanglements and ideological pressures interact to shape economic theory.

Sometimes this Soviet history reads as a description of a jigsaw puzzle that has all the essential elements but is never pieced together. In many cases, there were remarkable parallel developments suggesting similarities in the level of academic conversation in mathematics and mathematized disciplines both within and outside the Soviet Union. Nurtured by a system of specialized high school training and reinforced by the enthusiasm associated with cybernetics, Soviet mathematics formed a context in which a culture of formal argument could thrive in many applied fields, including, to some extent, economics. While the political regime changed over time, Soviet scholars found ways to navigate their institutional worlds while still producing useful–and sometimes pathbreaking–work. Moreover, despite many obstacles and delays, there was some communication between Soviet and non-Soviet research communities, especially with respect to the development of formal technical apparatuses. Soviet mathematicians and mathematical economists were aware of what was happening beyond the Iron Curtain, while their Western counterparts stayed informed about mathematical developments occuring in the USSR. This meant that some Soviet scholars were able to keep pace with Western academic research in economics, especially in its theoretical and formal areas. Many talented mathematicians were inspired by economic applications and enthusiastic about economic modeling. In the last decades of the USSR, an entire community of professional mathematical economists was formed.

---

[29] On this last episode see Komogorov's letter to Tikhomirov in Shiryaev (2006: 233). Girsanov's university lecture notes were published by Boris Polyak after Girsanov's untimely death, and quickly republished in English (Girsanov 1972). Girsanov's reviews in the journal of CEMI *Economics and Mathematical Methods* were also part of his pedagogical effort. These criticize the low quality of textbooks on the application of mathematics in the context of specific industries.

[30] On the various elements Kolmogorov, Dynkin, Skorokhod and others contributed to the foundations of mathematical finance, see Jarrow and Protter 2004.

And yet, even at this time, Soviet work in economics (Alexeev et al. 1992; Ericson 2019) never reached the level and prestige of Soviet mathematics. This paper has sought to narrate how and why this came about. In the decisive decades of the 1940s and 1950s, when economic theory embraced new mathematical formalisms (Debreu 1991), no joint conferences, no research visits, and no journal publications (Alexandrov 1996) were possible across the divide of the Iron Curtain. Bondareva's case vividly illustrated how the same problem affected the development of game theory in the 1970s and 1980s. The international mobility of most Soviet mathematicians and economists was extremely limited and the few rare exchanges remained largely inconsequential. The global adoption of their ideas has almost always followed the same pattern: formal technical contributions of Soviet mathematicians were applied to various fields of non-Soviet economic theory.

So strong mathematical research was not enough for economics to thrive. This story makes clear, however, that mathematics as context goes beyond a collection of technical results. Mathematicians were not "merely" developing tools for economists. As we have seen, von Neumann and Kantorovich were both deeply engaged in economics, and some of their mathematical work was directly inspired by economic applications. But as this article has argued, in the Soviet Union, similar avant-garde mathematics could not find proper economic applications. This discrepancy was clear already in 1970, when an *Econometrica* reviewer remarked, with some perplexity:

> The proper question is not why the country of Nemčinov, Kantorovič, and Novožilov has not produced 1000 x 1000 tableaux and Brookings SSRC models. It's just not yet that kind of country. The better question is perhaps why the country of Kolmogorov, Linnik, Dynkin, and many others has not produced something that elegantly supersedes the tenth Cowles Commission Monograph. (Fels 1970: 781)

It is telling, and was probably unknown by Fels at the time, that Kolmogorov and Linnik were active institutional supporters of mathematical economists. Nevertheless, their support did not outweigh the lack of a free and internationally open research environment.

A central element of this environment was official Soviet economic science, which was hostile to the direction taken by mainstream Western economics in the 1930s and 1940s. The deepest theoretical commitments of most Soviet economists led them to resist the micro-focused, behaviorally-oriented study of incentives and markets that increasingly defined Western economics at the time. In fact, in the USSR, the discipline never called itself *economics*: as far as economic theory was concerned, it remained *political economy* (of socialism). It thus faltered to accommodate Kantorovich's linear optimization and never embraced any game-theoretic ideas. Moreover, many of the influential Soviet economists played a role of academic watchdogs, helping to maintain the ideological purity of the discipline. (Self-)Censorship was widespread, enforced not only by powerful bureaucracies, but also by the Soviet economics profession.

The cases considered here draw these limitations into sharper focus. Because of endless ideological struggles, it even took Soviet linear programming, explicitly formulated as serving the needs of the planned economy, more than twenty years to take hold. When the Soviet research on optimal control and stochastic processes was at its zenith, it was largely ignored by economists in its country of provenance, but very quickly taken up in US macroeconomics. Likewise, Soviet game theorists did not believe their work had a particular relevance to economics. Unlike their Western colleagues, they did not have the chance to explore the possibilities engendered via interdisciplinary collaborations or potential applications in the social sciences.

Overall, the transfer of ideas between mathematics and economics happened very differently in the USSR and in the West. This article has argued that the different pathways of exchange and transmission can be explained by differing intellectual priorities and incentives, the minor role for economic theory in Soviet mathematical economists' training and research, and because of the obvious lack of space for application. The same mathematical results could be interpreted and applied differently. Context mattered.

The history of Soviet mathematicians and their engagements with economics helps us better understand how exactly context matters and how ideologically motivated limits to academic freedom and intellectual exchange can be seen – or, as Hicks would put it, photographically "developed" – as limits.

Literature

Alexandroff, Pavel S., and Heinz Hopf (1935). *Topologie. Bd I*. Berlin: Springer.

Alexandrov D. A. (1996). Why Soviet scientists stopped publishing abroad: How Soviet science became autonomous and isolated, 1914-1940, *Voprosy Istorii Esteststvoznaniya i Tekhniki* 3, 4–24. (In Russian.)

Alexeev, Michael, Clifford Gaddy, and J. Leitzel (1992). Economics in the Former Soviet Union. *Journal of Economic Perspectives*, 6(2): 137-148.

Aliprantis, Ch. D. and D. J. Brown (1983). Equilibria in markets with a Riesz space of commodities. *Journal of Mathematical Economics* 11, 189-207.

Aliprantis, Ch. and O. Burkinshaw (1978). *Locally Solid Riesz Spaces*. N.Y.: Academic Press.

Allen, R. G. D. (1950). The Work of Eugen Slutsky. *Econometrica*, *18*(3), 209–216.

Arrow, Kenneth J. and Frank H. Hahn (1971). *General Competitive Analysis*. San Francisco: Holden-Day.

Arrow, K. and M. Kurz (1970). *Public Investment, The Rate of Return, and Optimal Fiscal Policy*. Baltimore, MD: Johns Hopkins University Press.

Ashkenazi, V. (1961). *Applications of Game Theory to Military Affairs*. Moscow: Soviet Radio. (In Russian.)

Barnett, V. (2011). *E.E. Slutsky as Economist and Mathematician: Crossing the Limits of Knowledge*. L. and N.Y.: Routledge.

Bewley, T. (1972). Existence of equilibria in economies with infinitely many commodities. *Journal of Economic Theory* 4(3): 514-540.

Billingsley, P. (1968). *Convergence of Probability Measures*. N.Y. etc.: Wiley.

Bjork, S., Offer, A. & Söderberg, G. (2014). Time series citation data: the Nobel Prize in economics. *Scientometrics* 98, 185–196.

Bockman, J. and Bernstein, M. (2008). Scientific Community in a Divided World: Economists, Planning, and Research Priority during the Cold War. *Comparative Studies in Society and History*, 50(3): 581-613.

Boldyrev, I. (2020). Realities of Formalization: How Soviet Scholars Moved from Control Engineering to the General Theory of Choice, in: *Economics and Engineering: Institutions, Practices and Cultures, History of Political Economy*. Vol. 52, Annual Supplement (pp. 270–293), ed. by P. G. Duarte and Y. Giraud.

Boldyrev, I. (2024). The Frame for the Not-Yet-Existent: How American, European and Soviet Scholars Jointly Shaped Modern Mathematical Economics. *History of Political Economy, forthcoming*.

Boldyrev, I., Düppe, T. (2020). Programming the USSR: Leonid V. Kantorovich in context. *The British Journal for the History of Science*, 53(2), 255-278.

Boltyanski, V. (1999). The Maximum Principle – How it Came to Be? In: Boltyanski V., Martini, H., Soltan, V. Geometric Methods and Optimization Problems (pp. 204-230). Dordrecht: Kluwer.

Bondareva, O. N. (1962). Teoriia iadra v igre *n* lits (Theory of the core in the n-person game), (In Russian), *Vestnik LGU* 13, 141-142.

Bondareva, O. N. (1963/1968). "Nekotorye primeneniia metodor linejnogo programmirovaniia k teorii kooperativnykh igr" "Some applications of linear programming methods to the theory of cooperative games (In Russian)". *Problemy Kybernetiki*. 10: 119–139. Repr. in English in: *Selected Russian Papers in Game Theory*. Econometric Research Program, Princeton, 1968.

Bondareva, O. N. (1974). On Game-Theoretic Models in Economics. Leningrad: LGU. (In Russian.)

Borjas, G. J. and K. B. Doran (2012). The Collapse of the Soviet Union and the Productivity of American Mathematicians. *Quarterly Journal of Economics*, , 127 (3), 1143-1203.

Braumann, C. A. (2019). *Introduction to Stochastic Differential Equations with Applications to Modelling in Biology and Finance*. Hoboken, NJ: Wiley.

Brockwell, P. J. and R. A. Davis (1990). *Time Series: Theory and Methods*. N.Y.: Springer.

Burmeister, E. (2009). "Reflections." *History of Political Economy* 41, Suppl_1: 35–43.

Burmeister, E. and Dobell, A.R. (1970). *Mathematical Theories of Economic Growth*, New York: Macmillan.

Carvajalino J. (2021). Unlocking the Mystery of the Origins of John von Neumann's Growth Model. *History of Political Economy* 53 (4): 595–631.

Cass, D. (1965). Optimum Growth in an Aggregative Model of Capital Accumulation, *Review of Economic Studies* 32(3): 233–240.

Danilov, V. (1992). Implementation via Nash Equilibria. *Econometrica*, 60(1), 43–56.

Debreu, Gerard (1959). *Theory of Value: An Axiomatic Analysis of Economic Equilibrium*. New York: Wiley.

Debreu, G. (1986). Theoretic Models: Mathematical Form and Economic Content. *Econometrica*, 54(6): 1259–1270.

Debreu, G. (1991). The Mathematization of Economic Theory. *The American Economic Review*, 81(1), 1–7.

Demidov, Sergei S., and Boris V. Levshin (eds.). (2016). *The Case of Academician Nikolai Nikolaevich Luzin*, trans. Roger Cooke. Providence, RI: American Mathematical Society.

Dixit, A. and R. Pindyck (1994). *Investment under Uncertainty*. Princeton: Princeton University Press.

Dobrushin, R. L. 1968. "Description of a Random Field by Means of Conditional Probabilities and Conditions for Its Regularity." *Theory of Probability and Its Applications* 13: 197-224. (In Russian).

Dorfman, R. (1984). The discovery of linear programming. *Annals of the History of Computing* 6(3), pp. 283–295.

Dorfman, R., Samuelson, P. and R. Solow (1958). *Linear Programming and Economic Analysis*. N.Y. etc.: McGraw-Hill.

Düppe, T., & Weintraub, E. R. (2014). Siting the New Economic Science: The Cowles Commission's Activity Analysis Conference of June 1949. *Science in Context*, 27(3): 453-483.

Durlauf, S. N. (1993). Nonergodic Economic Growth. *Review of Economic Studies*, *60*(2), 349–366.

Ekeland, I. (2010). Notes on optimal transportation, *Economic Theory* 42: 437–459.

Erickson, Paul. (2015). *The World Game Theorists Made*. Chicago: The University of Chicago Press.

Ericson, R. (2019). The Growth and Marcescence of the "System for Optimal Functioning of the Economy" (SOFE). *History of Political Economy* 51 (S1): 155–179 (special issue: *Economic Knowledge in Socialism, 1945–89* ed. by Till Düppe and Ivan Boldyrev).

Faddeev, D. (1956). On the notion of entropy of finite probability distributions, *Uspekhi Mat. Nauk*, 11: pp. 227-231 (In Russian).

Fels, E. (1970). [Review of Mathematics and Computers in Soviet Economic Planning, by J. P. Hardt, M. Hoffenberg, N. Kaplan, & H. S. Levine]. *Econometrica*, 38(5), 779–781.

Föllmer, H. (1974). Random economies with many interacting agents, *Journal of Mathematical Economics*, 1(1): 51-62.

Fuchs, D. B. (2007). On Soviet Mathematics of 1950 and 1960s. In: Zdravkovska, Duren. 213-222.

Gale, D. (1960). *The Theory of Linear Economic Models*, McGraw-Hill book Co., New York. (Russ. translation 1963).

Galichon, A. (2016). *Optimal Transport Methods in Economics*. Princeton: Princeton University Press.

Galichon, A. (2021). The unreasonable effectiveness of optimal transport in economics. Prepared for the Proceedings of the 2020 World Congress of the Econometric Society. https://arxiv.org/pdf/2107.04700

Gamkrelidze, R.V. (1999). Discovery of the Maximum Principle. *Journal of Dynamical and Control Systems* **5**, 437–451.

Gardner, R. (1990). L.V. Kantorovich: the price implications of optimal planning. *Journal of Economic Literature* 28(2), pp. 638–648.

Gass, S. I. (1989). Comments on the history of linear programming. *IEEE Annals of the History of Computing* 11(2): 147-151.

Gerschenkron, A. (1978). Samuelson in Soviet Russia: A Report. *Journal of Economic Literature*, 16(2), 560–573.

Gikhman, I.I. and A. Skorokhod (1969). *Introduction to the Theory of Random Processes*. New York: Dover.

Giocoli, N. (2003). *Modeling Rational Agents: From Interwar Economics to Early Modern Game Theory*. Cheltenham: Edward Elgar.

Girsanov, I.V. (1960). On Transforming a Certain Class of Stochastic Processes by Absolutely Continuous Substitution of Measures. *Theory of Probability & Its Applications* 5(3): 285–301

Girsanov, I.V. (1972). *Lectures on Mathematical Theory of Extremum Problems*. Ed. by B.T.Polyak. Berlin etc.: Springer.

Gordon, L. A. (1992). *Dom (House)*. (In Russian). Saint-Petersburg: Neva.

Graham, Loren and Jean-Michel Kantor (2009). *Naming Infinity: A True Story of Religious Mysticism and Mathematical Creativity*. Cambridge, MA: Harvard University Press.

Hagemann, H, Kufenko, V, Raskov, D. (2016). Game theory modeling for the Cold War on both sides of the Iron Curtain. *History of the Human Sciences* 29(4-5): 99-124.

Hands, D. (1994). Restabilizing Dynamics: Construction and Constraint in the History of Walrasian Stability Theory. *Economics & Philosophy*, 10(2), 243-283.

Hansen, L.P., and T.J. Sargent (1980), "Formulating and estimating dynamic linear rational expectation models," *Journal of Economic Dynamics and Control* 2, 7-46.

Harrison, J. M., and D. M. Kreps (1979) Martingales and arbitrage in multiperiod securities markets, *Journal of Economic Theory* 20(3): 381-408.

Hicks, J. R. (1960). Linear Theory. *The Economic Journal*, 70(280): 671–709.

Hollings, C. (2013). 'The struggle against idealism: Soviet ideology and mathematics', *Notices of the American Mathematical Society* 60(11) (2013), 1448–1458

Jarrow, R. and P. Protter (2004). A short history of stochastic integration and mathematical finance: The early years, 1880–1970, A Festschrift for Herman Rubin Institute of *Lecture Notes-Monograph Series* 45: 75–91.

Kantorovich, L. V. (1937). Lineare halbgeordnete Räume, *Matematicheskiy sbornik* 2(1): 121–168.

Kantorovich, L. (1939/1960). *Mathematical methods in the organization and planning of production*. (In Russian). Leningrad: Leningrad University Press. Engl. translation in *Management Science* 6:366–422 (1959–60).

Kantorovich, L. V. (1942). "On the translocation of masses." *Doklady Akademii Nauk SSSR* 37(7–8):227–229. Transl. in *Management Science*, Vol. 5, No. 1 (Oct., 1958), pp. 1-4. Repr. in *Journal of Mathematical Sciences* 133(4):1381–1382 (2006).

Kantorovich, L. V. (1948). "On a problem of Monge." Uspekhi Mat. Nauk 3:225–226 (in Russian). English translation in *Journal of Mathematical Sciences* 133(4):1383 (2006).

Kantorovich, L. V. (1959/1965). Ekonomicheskii Raschet Nailuchshego Ispolzovania Resursov, Moscow: USSR Academy of Science (in Russian) (translated as *The Best Use of Economic Resources*, Cambridge, MA: Harvard University Press, 1965).

Kantorovich, L. V. (1966). Mathematical problems of optimal planning (In Russian). In: *Mathematical models and methods of optimal planning*. Novosibirsk: Nauka. (pp. 116–124).

Kantorovich, L. and Rubinstein, G. (1958). "On a space of completely additive functions." (In Russian). Vestn. Leningrad Univ. 13(7):52–59.

Kantorovich, Vsevolod L., Kutateladze, Semion S. and Fet, Yakov I. (eds.) (2004). Leonid V. Kantorovich: Chelovek i uchenyii (Leonid V. Kantorovich: Man and Scientist), vol. 2, Novosibirsk: Scientific Publishing Center of RAS.

Karlin, Samuel (1959/1964). *Mathematical Methods and Theory in Games, Programming, and Economics*. Reading, MA: Addison-Wesley. (Russ. Translation in 1964).

Khinchin, A. Ya. (1933). Über stationäre Reihen zufälliger Variablen'. *Matematicheskii Sbornik*
40, 124–128.

Khinchin, A.Ya. (1934). 'Korrelationstheorie der stationären stochastischen Prozesse'. *Mathematische Annalen* 109, 604–615.

Khinchin, A. Ya. (1953). *The entropy concept in probability theory*. Uspekhi Mat. Nauk, 8: pp. 3-20. (In Russian.)

Khinchin, A. Ya. (1957). *Mathematical Foundations of Information Theory*. N.Y.: Dover.

Khintchine, A. Y (1932). "Mathematical theory of a stationary queue". *Matematicheskii Sbornik*. 39 (4): 73–84. (In Russian.)

Khintchine, A. Ya. (1960). *Mathematical Methods in the Theory of Queuing*. London: Charles Griffin & Co.

Kjeldsen, T. H. (2001) John von Neumann's Conception of the Minimax Theorem: A Journey Through Different Mathematical Contexts. *Archive for History of Exact Sciences* 56, 39–68.

Klein J. L. (1999). The Rise of "Non-October" Econometrics: Kondratiev and Slutsky at the Moscow Conjuncture Institute. *History of Political Economy* 31 (1): 137–168.

Kolmogorov, A. N. (1931a). Sur le problème d'attente, *Matematicheskii Sbornik* 38 (1—2), 101—106.

Kolmogorov, Andrei (1931b). "Über die analytischen Methoden in der Wahrscheinlichkeitsrechnung" *Mathematische Annalen* **104**: 415–458.

Kolmogorov [Kolmogoroff], A.N. (1933). *Grundbegriffe der Wahrscheinlichkeitsrechnung*. Berlin: Julius Springer.

Kolmogorov, A. N. (1941a). Stationary sequences in Hilbert space. *Herald of Moscow State University*. 2 (6) 1-40. (In Russian.)

Kolmogorov, A.N. (1941b). Interpolation and extrapolation of stationary random sequences. *Bulletin de l'Academie des Sciences de U.S.S.R., Ser. Mathematics* 5: 3–14. (In Russian.)

Kolmogorov, A.N. (1948/2002), E.E. Slutsky. An obituary. *Math. Scientist*, vol. 27, 2002, pp. 67 – 74.

Koopmans, T. C. (1951) (ed.): *Activity Analysis of Production and Allocation*. Cowles Commission Monograph, 13, New York:Wiley.

Koopmans, Tjalling (1957). *Three Essays on the State of Economic Science*. New York: McGraw-Hill.

Kornai, János (1979). Resource-Constrained versus Demand-Constrained Systems. *Econometrica*, *47*(4), 801–819.

Kornai, János (1980). *Economics of Shortage*, Amsterdam: North Holland

Kornai, J., & Th. Lipták. (1965). Two-Level Planning. *Econometrica*, *33*(1), 141–169.

Kreps, D. M. (1981). Arbitrage and equilibrium in economies with infinitely many commodities, *Journal of Mathematical Economics*, 8(1): 15-35.

Kulakovskaja T. E., Naumova N. I. (1992). Olga Nikolajevna Bondareva. 1937-1991, *International Journal of Game Theory*. 20(4):309—312.

Magill, M. J. P. (1970). *On a General Economic Theory of Motion*, N.Y.: Springer-Verlag.

Maschler, M. (1992). The Bargaining Set, Kernel, and Nucleous. In: *Handbook of Game Theory*, Volume 1, Edited by R.J Aumann and S. Hart. N.Y.: Elsevier.

Massey, F. J. (1951). The Kolmogorov–Smirnov test of goodness of fit. *Journal of the American Statistical Association* 46, 68–78.

Mills, T. C. (2011). *The Foundations of Modern Time Series Analysis*. N.Y.: Palgrave Macmillan.

Minkowski, H. (1911). *Gesammelte Abhandlungen*. Bd. 2. Leipzig, Berlin: B. G. Teubner.

Mirowski, Philip (2002). *Machine Dreams: Economics Becomes a Cyborg Science*. New York: Cambridge University Press.

Mises, L. von (1920). Die Wirtschaftsrechnung im sozialistischen Gemeinwesen. *Archiv für Sozialwissenschaft und Sozialpolitik* 47, 86-121.

Mityagin, B. S. (1972). "Zametki po matematicheskoi ekonomike", *Uspekhi matematicheskikh nauk*, 27:3(165): 3–19; *Russian Math. Surveys*, 27:3, 1–19

Morgenstern, O. (1976). The Collaboration Between Oskar Morgenstern and John von Neumann on the Theory of Games. *Journal of Economic Literature 14*(3): 805-816.

Myerson, R. (1992). In Memoriam: Olga Bondareva, *Games and Economic Behavior* 4: 324.

Myerson, R. B. (1999). Nash Equilibrium and the History of Economic Theory. *Journal of Economic Literature*, 37 (3): 1067-1082.

Pesch, H. J. and Plail, M. (2009). The Maximum Principle of Optimal Control: A History of Ingenious Ideas and Missed Opportunities. *Control and Cybernetics* 38, No. 4A, 973-995.

Phillips, P. C. B. & Perron, P. (1988). Testing for a Unit Root in Time Series Regression. *Biometrika*, *75*(2), 335–346.

Polyak, B. T. (2002). History of Mathematical Programming in the USSR: Analyzing the Phenomenon. *Math. Program.* 91, 401–416.

Pontryagin, L. S. (1998). Biography of Lev Semyonovich Pontryagin, mathematician, compiled by himself. Moscow: Prima. (In Russian).

Pontryagin, L.S., Boltyanskii, V. G., Gamkrelidze, R. V. , Mishechenko, E. F. (1962). *The Mathematical Theory of Optimal Processes*. New York/London: Wiley.

Raghavan, T. E. S. (2002). Non-Zero-Sum Two-Person Games In: *Handbook of Game Theory*, Volume 3, Edited by R.J Aumann and S. Hart. N.Y.: Elsevier.

Ramsey, Frank P. (1928). A Mathematical Theory of Saving, *The Economic Journal*, 38(152): 543–559.

Robinson, T. W. (1970). Game Theory and Politics: Recent Soviet Views. *Studies in Soviet Thought*, 10(4), 291–315.

Romanovskii, I.V. (1962), "Reduction of a game with complete memory to a matrix game", *Soviet Mathematics* 3:678-681.

Roth, A. E. and I. Erev (1995). Learning in Extensive-form Games: Experimental Data and Simple Dynamic Models in the Intermediate Term, *Games and Economic Behavior* 8(1): 164-212.

Rozanov, Yu.A. (1967). *Stationary Random Processes*. San Francisco: Holden-Day.

Rubinstein, A. (1986). Finite Automata Play the Repeated Prisoner's Dilemma, *Journal of Economic Theory* 39(1): 83-96.

Samuelson, P. (1947). *Foundations of Economic Analysis*. Cambridge, MA: Harvard University Press.

Samuelson, P. A. (1972). Maximum Principles in Analytical Economics. *American Economic Review*, *62*(3), 249–262.

Samuelson, P.A. and Solow, R.S. (1956). "A complete capital model involving heterogeneous capital goods", *Quarterly Journal of Economics* 70:537–62.

Schwartz, B.L. (1989). The invention of linear programming. *IEEE Annals of the History of Computing* 11(2): 145-147.

Seneta, E. (2004). Mathematics, religion, and Marxism in the Soviet Union in the 1930s, *Historia Mathematica*, Volume 31, Issue 3: 337-367.

Sent, Esther-Mirjam. 1998. *The Evolving Rationality of Rational Expectations: An Assessment of Thomas Sargent's Achievements*. New York: Cambridge University Press.

Shafer, G. and Vladimir Vovk (2006). The Sources of Kolmogorov's "Grundbegriffe" *Statistical Science*, Vol. 21, No. 1, pp. 70-98.

Shapley, Lloyd. 1967. On balanced sets and cores. *Navel Research Logistics Quarterly* 9, 45–48.

Shapley, Lloyd S. and Martin Shubik (1954). A method for evaluating the distribution of power in a committee system. *American Political Science Review* 48: 787–792.

Shell, Karl. (1987) "Hamiltonians" in The New Palgrave: A Dictionary of Economics (J. Eatwell, M. Milgate and P. Newman, eds.), Vol. 2, New York: Macmillan, 1987, 588-590.

Shiryaev, A. (2006). Kolmogorov in the Recollections of His Students. Moscow: Moscow Center for Continuous Mathematical Education. (In Russian.)

Skorokhod, A. V. (1965). *Studies in the Theory of Random Processes*. Reading, MA: Addison-Wesley.

Slutsky, E.E. (1915/1952). On the Theory of the Budget of the Consumer. In A.E.A. Readings in Price Theory, ed. by G. J. Stigler and K. E. Boulding (pp. 27-56) Chicago: Richard D. Irwin. Originally published as "Sulla teoria del bilancio del consumatore" in *Giornale degli Economisti* 51 : 1-26.

Slutzky, E. (1937). The Summation of Random Causes as the Source of Cyclic Processes. *Econometrica*, 5(2), 105–146.

Stengel, B. von (2002). Computing Equilibria for Two-Person Games. In: *Handbook of Game Theory*, Volume 3, Edited by R.J Aumann and S. Hart. N.Y.: Elsevier.

Stokey, N. (2009). *The Economics of Inaction: Stochastic Control Models with Fixed Costs*. Princeton: Princeton University Press.

Tsetlin, M.L. (1973). *Automaton Theory and Modeling of Biological Systems*. New York and London: Academic Press.

Weintraub, E. R. (ed.) (1992). *Toward a History of Game Theory*. Durham, NC and L.: Duke University Press.

Uzawa, H. (1965). Optimum Technical Change in An Aggregative Model of Economic Growth. *International Economic Review*, *6*(1), 18–31.

Vasil'ev, V.A. and E. B. Yanovskaya (2006). Introduction. In: Driessen T. H. et al. (eds.). *Russian Contributions to Game Theory and Equilibrium Theory*. Berlin and N. Y.: Springer.

Vershik, A.M. (2013). Long History of the Monge-Kantorovich Transportation Problem. *Mathematical Intelligencer* 35, 1–9.

Volkonskii, V. (1967). *The Model of Optimal Planning and the Interrelations of Economic Parameters*. (In Russian). Moscow: Nauka.

von Neumann, John (1928). Zur Theorie der Gesellschaftsspiele. *Mathematische Annalen* 100: 295–320.

von Neumann, J. (1937/1945). Über ein ökonomisches Gleichungssystem und eine Verallgemeinerung des Brouwerschen Fixpunktsatzes, in K. Menger, ed. *Ergebnisse eines mathematischen Kolloquiums, 1935-36*. Engl. Translation: A Model of General Economic Equilibrium. *Review of Economic Studies* 13 (1): 1-9.

Vorob'ev, N.N. (1958), "Equilibrium points in bimatrix games", *Theory of Probability and its Applications* 3:297-309.

Vorob'ev, N. N. (1977). *Game Theory*. Transl. and supplemented by S. Kotz. New York and Berlin: Springer.

Whittle, P. (1963). *Prediction and Regulation*. London and Princeton: English University Press and Van Nostrand.

Wald, A. (1949) 'Statistical decision functions', *Annals of Mathematical Statistics*, 20:165-205.

Whitney, H. (1988). Moscow 1935: Topology Moving Toward America. In: Peter Duren (ed.) A century of Mathematics in America. Vol 1. Providence, RI: AMS. (pp. 97-117).

Wiener, N. (1949). *Extrapolation, Interpolation and Smoothing of Stationary Time Series, with Engineering Applications*. New York: Wiley.

Wold, H. (1938). *A Study in the Analysis of Stationary Time Series*. Stockholm: Almqvist & Wiksell.

Wooders, M. (2008). "Bondareva, Olga (1937–1991)." *The New Palgrave Dictionary of Economics*. 2nd Edition. Eds. S. N. Durlauf and L. E. Blume. Palgrave Macmillan.

Wulwick, N. J. (1995). The Hamiltonian formalism and optimal growth theory. In: Ingrid H. Rima (ed.) *Measurement, Quantification and Economic Analysis: Numeracy in Economics*. L. and N.Y.: Routledge.

Yaglom, A. M. (1962). *An Introduction to the Theory of Stationary Random Functions*. N.Y.: Prentice-Hall.

Yasny, N. (1972). *Soviet Economists of the Twenties. Names to be Remembered*. Cambridge: Cambridge University Press.

Yurchak, Alexei. (2005). *Everything Was Forever, Until It Was No More: The Last Soviet Generation*. Princeton: Princeton Unviersity Press.

Zauberman, A. (1975a). *The Mathematical Revolution in Soviet Economics*, Oxford etc.: Oxford University Press.

Zauberman, A. (1975b). *Differential Games and Other Game-Theoretic Topics in Soviet Literature: A Survey*. New York: New York University Press.

Zauberman, A. (1976). *Mathematical Theory in Soviet Planning*. London etc.: Oxford University Press.

Zdravkovska, S. and P. L. Duren (eds.) (2007). *Golden Years of Moscow Mathematics*. Providence, R.I.: American Mathematical Society and London Mathematical Society, 2nd ed.

Appendix

**Soviet Mathematical Economics: A Very Short Introduction (With Some Illustrations)**

This Appendix takes a closer look at Soviet mathematical economics (rather than mathematics; see also Introduction). It elaborates on its contexts and discusses in more detail some exemplary contributions, focusing predominantly on the postwar years.

The term "mathematical economics" is not very specific as most postwar economics was mathematized. Perhaps one could say that some areas of economic theory are primarily motivated by the interest to solve a mathematical problem, and as such could be called "mathematical economics." In the Soviet context, however, the term does have a more specific meaning. It captures any work that was informed by mathematical techniques and thus deviated from what was largely taught and researched in the economic departments within the system of higher education or the Soviet Academy of Sciences.

Joseph Stalin was responsible not only for the repression of many economists, including those associated with the Conjuncture Institute (see Introduction), but also for instilling an overarching mistrust toward quantitative reasoning in the social sciences. For decades, this sort of research was ideologically suspect. Only after Stalin's death in 1953, with the softening of political regime (the epoch of the "Thaw"), could the discussion around applying mathematical methods in economics be re-opened (Ellman 1973; Gerovitch 2002; Leeds 2016; Boldyrev and Kirtchik 2017). The first study program was organized in 1959 at the economics department of Leningrad University, where Kantorovich was active. In 1960, a similar program opened at Moscow University, initiated by Nemchinov. In 1963, the Central Economic-Mathematical Institute (CEMI) was created. It became the major mathematical economics research institution in the Soviet Union. In 1964, with the help of Girsanov (see Section 4.2) and others, a separate division of economic cybernetics was created there, subsequently emulated by other universities across the whole country. The year 1965 saw the launch of the first professional journal, *Economics and Mathematical Methods*, institutionalizing the field. Formally, this research was relegated to the separate disciplines like "mathematical methods" and "economic cybernetics," with no claims to become "economic theory." This defined the subaltern status of the field.

Note that most of Soviet science was being done not at the universities (although many scholars had regular teaching responsibilities) but in the research centers associated with the Academy of Sciences. Typically, the research was published in the home journals of the relevant institutes. Some of it was translated, but not widely read. Soviet scholars almost never submitted their work to international economic journals.

There were two important social tendencies that defined the context for mathematical economics in the USSR. First, quite generally, mathematics and other less "controversial" technical disciplines attracted talented people (whatever definition of talent one adopts), who otherwise would have pursued a career in the humanities or in the social sciences. The fields were seen as too ideological and did not facilitate open and critical conversation. Mathematical economics, on the contrary, provided a space that allowed for some more freedom to discuss alternative forms of economic organization. The second tendency was the infamous academic antisemitism, for which, sadly, various mathematicians, including Pontryagin, bear some responsibility. From the mid-1960s to the mid-1980s, an academic career in mathematics, especially at the major centers, like Steklov Mathematical Institute and the *mekhmat* of Moscow University, became difficult for Jews. They were discriminated against and had lots of troubles in getting into a PhD programme, receiving academic appointments, and traveling to international conferences. While the antisemitic discrimination has been well documented (see the materials in Shifman 2005), there does not seem to be any systematic studies of the general tendency to avoid humanities and social sciences. Nonetheless, it repeatedly pops up in the memoirs of mathematicians.

For example, consider Anatole Katok, a mathematician of a Jewish origin, who was, before his emigration in 1978, a researcher at CEMI. His biography illustrates *both* tendencies. Katok

> explained that his choice of mathematics as a vocation was influenced by the relative freedom mathematicians enjoyed because their discipline was least affected and controlled by ideological impositions [...] From the late 1960s antisemitism and suppression of liberal thought grew at Moscow State University, and almost no Jews were accepted as students or faculty. So, Katok instead assumed an appointment at the [...] CEMI [...], which allowed him to combine work on mathematical problems in economics, if any, with research in pure mathematics (Hasselblatt 2019: 711).

The new community of mathematical economists integrated scholars of the older generation. Thus, the former members of the Conjuncture Institute who survived the repressions were able to return to research. These included the statistician Yakov Gerchuk, or the statistician and cliometrician Albert Vainshtein, as well as other mathematically minded economists of the previous generations, such as Alexander Lur'e and Viktor Novozhilov. All these individuals played some role in creating a research environment that should have been different from the vacuum of the previous decades.

The work of Alexander Konüs (1895-1990), another surviving economist, is of particular interest to illustrate how huge this vacuum was. Konüs became famous for his early study on the cost-of-living index (Konüs 1924/1939).[31] His work was probably the most fruitful collaboration of an economist and a mathematician in the Soviet Union between the wars. Konüs worked with the Moscow University geometer Sergei Byushgens. Their paper (Konüs and Byushgens 1926) was rediscovered in the 1970s by W. E. Diewert (1976), a Canadian economist of Russian origin, who studied the language and could read their text in Russian. It has been recently translated into English. In the preface to the translation, Diewert and Zelenyuk (2023) write that the paper contained several important theoretical results. For example, the early demonstration of microeconomic duality theory makes the paper not only "a landmark in the history of index number theory," but also an important document in the history of microeconomics.

CEMI enjoyed relative intellectual autonomy, providing a home for many scholars who would otherwise have had troubles getting academic jobs. Its director Nikolai Fedorenko promoted CEMI as an important organizational project and was for years defending, not without success, both the institute and the specific academic program CEMI suggested for policymakers. This program came to be associated with the so-called SOFE (system of optimal functioning of the economy), a multi-level system of models that was supposed to inform planning.[32] Although the relations between CEMI and Gosplan, the Soviet Planning Committee, were strained, the two organizations did manage to collaborate, and CEMI proposals informed some planning initiatives. Perhaps the most important bridging project was that of ASPR (*Automatized System of Planning Calculations*), which sought to integrate insights from input-output and optimization models into the planning process (Urinson 1986). During the *perestroika* (1985-1991) and the first years of economic reforms, or, as Makarov (1988: 459) put in at the 100th AEA meeting, during "revolutionary changes [,] a transition from an excessively stable and rigid economic system to one which is much more flexible," the economists affiliated with or coming from CEMI turned out to be influential economic experts and policymakers.

Theoretical and applied work on economic modeling was also done at the Gosplan research institute, at the Chief Computer Center of the Academy of Sciences,[33] and at the Institute of Control Problems. Important institutions outside Moscow were located in Novosibirsk.[34] Kantorovich was working there in the 1960s, but after he left for Moscow, many mathematicians and mathematical economists, including Vladimir Bulavskii, Valeri Marakulin, Leonid Polishchuk, Alexander Rubinov, Gennady Rubinstein, and Valeri Vasil'ev, continued this research. In Kyiv, the Institute of Cybernetics of the Ukrainian Academy of Sciences was another important center.

---

[31] The paper was known to econometricians of that time via references, but not in full. It was Henry Schultz, American demand theorist and statistician, who, in the mid 1930s, organized the translation of the paper he could not read, used it in his teaching, and suggested publishing it in *Econometrica*. (Schultz also initiated the publication of Slutzky (1937)).

[32] The first formulations were given by Volkonskii (1967) and Katseneligenboigen et al. (1969), see the history in Ericson 2019.

[33] There, the department of mathematical economics was founded in 1968 by Nikita Moiseev. On Moiseev see Rindzevičiūtė (2016).

Although international collaborations with Soviet economists were very limited, there were some exceptions. For example, Martin Weitzman co-authored papers with Soviet authors, frequently visited the USSR, and generally belonged to those few American mainstream economists who managed to maintain long-term contacts with the Soviet mathematical economics community.[35] Other major economists belonging to this group were Tjalling Koopmans (Düppe 2016) and David Gale.

In the beginning of the 1980s, Yuri Yermoliev, an operations researcher and expert in stochastic programming from the Kyiv Institute of Cybernetics, and another Soviet Ukrainian mathematician, Yuri Kaniovski, collaborated with W. Brian Arthur at the International Institute for Applied Systems Analysis (IIASA) in Austria, while Arthur was doing his now famous research on increasing returns (published later as Arthur 1989). This collaboration–initiated at the institution created precisely to further academic exchange across the Iron Curtain (Rindzevičiūtė 2016)–illustrates the potential carried by the international division of academic labor between Western economics and Soviet mathematics.

What kind of research prevailed in Soviet mathematical economics? One cannot possibly do justice to all the relevant developments over several decades. However, a bird's-eye view can provide a general idea of what the field was about.

It would be fair to say that a lot of intellectual effort, beginning from the end of the 1950s, was directed at input-output modeling. This technique was connected both with the idea of improving economic calculation, considered as a basic rationale behind the "mathematical methods," and with the promise of practical application. Although the practical results were reported to be quite modest (Tretyakova and Birman 1976), input-output techniques became the standard element in the university curricula in "economic cybernetics" and in applied research. Along with input-output, many mathematical economists were busy developing the efficiency criteria for investment–something not easy in the context where openly discussing capital theory was risky. Another CEMI economist, Boris Mikhalevski, explored repressed inflation. But studying the actual state of the Soviet economy was also subject to censorship.

At CEMI, apart from the work that promised a practical import, a lot of attention was paid to improving optimization techniques (see the overview in Polyak 2002). CEMI served as a temporary home for some bright mathematicians, including Eugene Dynkin, Gennadi Henkin, Anatole Katok, and Boris Mityagin.[36] Some of them influenced this research agenda and contributed to defining standards of rigor in theoretical modeling. Overall, Soviet operations researchers made significant contributions, as evidenced by the work of David Yudin, Arkadi Nemirowski, Yuri Nesterov, and Leonid Khachiyan, who received multiple international awards.[37]

---

[34] The Institute of Economics and Organization of Industrial Production and the Mathematical Economics division at the Institute of Mathematics were both of the Siberian Division of the Academy of Sciences.

[35] For example, Kantorovich's (1965) book was used in Weitzman's 1974 MIT course on core micro theory https://www.irwincollier.com/m-i-t-core-micro-theory-resource-allocation-price-system-weitzman-1974/

[36] On Dynkin's contributions and on the links with further literature, see Evstigneev (2000).

In the beginning of the 1960s, Kantorovich's student in Novosibirsk, Valery Makarov, became interested in optimal growth theory and general equilibrium analysis (Makarov and Rubinov 1973/1977, see the review by Gale 1978). Makarov moved to Moscow in 1985 to become the director of CEMI. The Novosibirsk scholars did mathematically sophisticated work on optimization, general equilibrium, and cooperative games. At CEMI, Arkin and Evstigneev (1987) reformulated optimal growth theory in stochastic terms, a contribution that played some role in this literature (see Brock and Dechert 2010 for an overview).

Victor Polterovich, another CEMI economist, was particularly active as a general equilibrium theorist. He formulated a criterion for the monotonicity of aggregate demand functions (Mitiushin and Polterovich 1978); provided a version of a synthesis between optimal growth and general equilibrium theory (Polterovich 1983); and used disequilibrium modeling to understand the realities of shortage (Polterovich 1990; 1993). In fact, disequilibrium analysis in the USSR was pioneered by Emmanuil Braverman (1972), who was working at the Institute of Control Problems. Another influential group at the same institute, led by Mark Aizerman, engaged, since the 1970s, in research on voting and abstract choice theory.[38] The Aizerman group was, perhaps, one of the most internationalized research collectives in Soviet mathematical economics, hosting regular seminars and being visited by major social choice theorists and game theorists of the time.

When asked in the 1980s, some of these protagonists said that at the time, mathematics was perceived as the way to improve the general culture of economic planning (Katseneligenboigen 1981) and to infuse Soviet economics with some rationality (Volkonskii 1989). Indeed, mathematical formalisms, once accepted, helped further the discussion – at least academically. Overall, the work of Soviet mathematical economists remained quite abstract and remote from actual applications.[39] But sometimes it was brilliant work, still present in today's formal economic theory. In this short overview, I have tried to demonstrate that exploring the intellectual legacy of this diverse field alongside the experiences of Soviet mathematical economists and the context of their work is an instructive endeavor in the history of economic ideas.

Arkin, V. & I. Egstineev (1987). *Stochastic Models of Control and Economic Dynamics*. New York: Academic Press

---

[37] In 1982, Yudin, Nemiroski and Khachiyan received Fulkerson Prize in discrete mathematics, Nemirovski and Nesterov received Dantzig Prize in mathematical programming (in 1991 and 2000, respectively), and the John von Neumann Prize (in 2003 and 2009), the major award in operations research.

[38] On the career of Polterovich and his work in demand theory and disequilibrium modeling, see Boldyrev and Kirtchik (2014), Boldyrev (2023); Braverman's work is covered in Kirtchik (2019); on Aizerman and his group (in particular, Andrei Malishevski and Fuad Aleskerov) see Boldyrev (2020).

[39] Yuri Kabanov, currently a scholar in mathematical finance (see also 4.2), worked in the CEMI lab headed by Vadim Arkin and focused on stochastic calculus. In a characteristic episode, he recalls having encountered a paper by Merton (1975). 'I asked Vadim [Arkin] […]: "What is a portfolio?" He said: "Yura, forget about it. […] In capitalist countries, there are portfolios, in socialist economies, it is not an important issue, you cannot be promoted at our institute [if you focus on that] […] one should do something more applied, more economic"' (Kabanov 2010). Indeed, in a country without financial markets, portfolio theory was entirely impractical.

Arthur, W. Brian. (1989). Competing Technologies, Increasing Returns, and Lock-In by Historical Events. *The Economic Journal*, *99*(394), 116–131.

Boldyrev, I. (2023). Victor Polterovich. In: Vernengo, M., Caldentey, E.P., Rosser Jr, B.J. (eds) *The New Palgrave Dictionary of Economics*. Palgrave Macmillan, London. https://doi.org/10.1057/978-1-349-95121-5_3142-1

Boldyrev, I. Kirtchik. O. (2014). General Equilibrium Theory Behind the Iron Curtain: The Case of Victor Polterovich, *History of Political Economy*. 46(3): 435-461.

Boldyrev, I., and Kirtchik, O. (2017). "The Cultures of Mathematical Economics in the Postwar Soviet Union: More than a Method, Less than a Discipline." *Studies in History and Philosophy of Science Part A* 63 (6): 1–10.

Braverman, E. (1972). Model proizvodstva s neravnovesnymi tsenami. [Model of Production with Disequilibrium Prices], *Ekonomika i matematicheskie metody*, 8, 40-64 (in Russian).

Brock, W.A., Dechert, W.D. (2010). growth models, multisector. In: Durlauf, S.N., Blume, L.E. (eds) Economic Growth. The New Palgrave Economics Collection. Palgrave Macmillan, London. (pp. 127 -132).

Diewert, W.E. (1976). Exact and superlative index numbers, *Journal of Econometrics*, 4(2): 115-145.

Diewert, E., Zelenyuk, V. (2023). On the problem of the purchasing power of money by A. A. Konüs and S. S. Byushgens: translation and commentary. *Journal of Productivity Analysis*. https://doi.org/10.1007/s11123-023-00696-x

Düppe, T. (2016). Koopmans in the Soviet Union: a travel report of the summer of 1965. *Journal of the History of Economic Thought* 38(1): 81–104.

Ellman, M. 1973. *Planning problems in the USSR: The contribution of mathematical economics to their solution 1960-1971*. Cambridge: Cambridge University Press.

Evstigneev, I. V. (2000). Dynkin's work in mathematical economics. In: *Selected papers of E. B. Dynkin with Commentary*. A.A. Yushkevich, G. M. Seitz, A.L.Onishchik (eds.). Providence, RI: AMS, and Cambridge, MA: International Press. pp. 793-795.

Gale, David (1978). Review of: Mathematical theory of economic dynamics and equilibria, by V. L. Makarov and A. M. Rubinov. *Bulletin of the American Mathematical Society*, 84, 665-671.

Gerovitch, S. (2002). *From Newspeak to Cyberspeak: A History of Soviet Cybernetics*. Cambridge, MA: MIT Press.

Hasselblatt, Boris (2019). Anatole Katok - A Half-Century of Dynamics, *Notices of the AMS* 66(5): 708 – 719.

Kabanov, Y. (2010). Mathematical finance and mathematics from finance. A talk at the International Symposium *Visions in Stochastics (Leaders and their Pupils)* Moscow, November 1, 2010. https://www.mathnet.ru/php/presentation.phtml?option_lang=eng&presentid=2554

Katseneligenboigen, Aron (1981). Interview conducted by Eugene Dynkin with Aron Katsenelinboigen on July 5, 1981. https://hdl.handle.net/1813/17453

Katsenelinboigen, A., Lakhman, I., & Ovsienko, Yu. (1969). Optimalnost' i tovarno-denezhnye otnosheniia [Optimality and Commodity-Money Relations]. Moscow: Nauka (in Russian).

Kirtchik, O. (2019). From Pattern Recognition to Economic Disequilibrium: Emmanuil Braverman's Theory of Control of the Soviet Economy. *History of Political Economy* 51 (S1): 180–203.

Konüs, A.A. (1924). "The Problem of the True Index of the Cost of Living," The Economic Bulletin of the Institute of Economic Conjuncture, Moscow, No. 9-10 (36-37), pp. 64-71, translated in *Econometrica* 7, (1939), 10-29.

Konüs, A.A., and S.S. Byushgens (1926) `K probleme pokupatelnoi cili deneg' (On the problem of the purchasing power of money'), Voprosi Konyunkturi II(1) (supplement to the Economic Bulletin of the Conjuncture Institute), 151-72.

Leeds, A. (2016). Dreams in Cybernetic Fugue: Cold War Technoscience, the Intelligentsia, and the Birth of Soviet Mathematical Economics. *Historical Studies in the Natural Sciences* 46 (5): 633–668.

Makarov, V. L. (1988). On the Strategy for Implementing Economic Reform in the USSR. *The American Economic Review*, *78*(2), 457–460.

Makarov, V., & Rubinov, A. (1977). *Mathematical theory of economic dynamics and equilibria* (Translation of *Matematicheskaya teoriya ekonomicheskoi dinamiki i ravnovesiya*. Moscow: Nauka, 1973). New York, Heidelberg, Berlin: Springer.

Merton, R. C. (1975). Optimum Consumption and Portfolio Rules in a Continuous-Time Model, In: W.T. Ziemba, R.G. Vickson (eds.), *Stochastic Optimization Models in Finance* (pp. 621-661). N.Y.: Academic Press.

Mitiushin, L. and V. Polterovich (1978) Criterion for monotonicity of demand functions. *Ekonomika i matematicheskie metody [Economics and Mathematical Methods]* 14(1): 122–128. (In Russian).

Polterovich, V. M. (1983). Equilibrium Trajectories of Economic Growth. *Econometrica*, *51*(3), 693–729.

Polterovich, V. M. (1990). *Economic equilibrium and economic mechanism*. Moscow: Nauka. (In Russian).

Polterovich, V.M. 1993. Rationing, queues, and black markets. *Econometrica* 61 (1): 1–28.

Rindzevičiūtė, Egle (2016). *The Power of Systems: How Policy Sciences Opened Up the Cold War World*. Ithaca and London: Cornell University Press.

Shifman, M. (ed.) (2005). *You Failed Your Math Test, Comrade Einstein*. Singapore: World Scientific.

Tretyakova A. and I. Birman (1976). "Input-output Analysis in the USSR." *Soviet Studies* 28, no. 2: 157-186.

Urinson, Y. M. (1986). *Perfecting the Technology of National Economic Planning* (In Russian.) Moscow: Ekonomika.

Volkonskii, V. (1989) Interview conducted by Eugene Dynkin with Viktor Volkonsky, September 16, 1989. http://dynkincollection.library.cornell.edu
https://ecommons.cornell.edu/items/f90fdfec-602a-47f1-a887-9630ef5baf9b